\documentclass[11pt,twoside]{amsart}
\usepackage{amsmath, amsthm, amscd, amsfonts, amssymb, graphicx, color}
\usepackage[bookmarksnumbered, plainpages]{hyperref}
\usepackage{enumitem}
\numberwithin{equation}{section}
\def\pn{\par\noindent}

\newtheorem{thm}{Theorem}[section]

\newtheorem{lem}[thm]{Lemma}
\newtheorem{cor}[thm]{Corollary}
\newtheorem{prop}[thm]{Proposition}
\newtheorem{rmk}[thm]{Remark}
\newtheorem{exmp}[thm]{Example}
\newtheorem{problem}[thm]{Problem}

\theoremstyle{remark}

\begin{document}
%------------------------------------------------------------------------------------%
%%Don not change any thing in this part
\hskip -0.2 cm
\begin{tabular}{c r}
\vspace{-0.6cm}
\href{http://www.theoryofgroups.ir}{\scriptsize  \rm www.theoryofgroups.ir}\\
\end{tabular}
\hskip 2 cm
\begin{tabular}{l}
\hline
\vspace{-0.2cm}
\scriptsize \rm\bf International Journal of Group Theory\\
\vspace{-0.2cm}
\scriptsize \rm ISSN (print): 2251-7650, ISSN (on-line): 2251-7669 \\
\vspace{-0.2cm}
\scriptsize Vol. {\bf\rm x} No. x {\rm(}201x{\rm)}, pp. xx-xx.\\
\scriptsize $\copyright$ 201x University of Isfahan\\
\hline
\end{tabular}
\hskip 2 cm
\begin{tabular}{c c}
\vspace{-0.1cm}
\href{http://www.ui.ac.ir}{\scriptsize \rm www.ui.ac.ir}\\
\vspace{-1cm}
\end{tabular}
\vspace{1.3 cm}

%------------------------------------------------------------------------------------%

\title{Some problems about products of conjugacy classes in finite groups}
\author{Antonio Beltr\'an, Mar\'ia Jos\'e Felipe and Carmen Melchor*}

\thanks{{\scriptsize
\hskip -0.4 true cm MSC(2010): Primary: 20E45; Secondary: 20C15.
Received: dd mmmm yyyy, Accepted: dd mmmm yyyy.\\
$*$Corresponding author}}
\maketitle

%------------------------------------------------------------------------------------%
%This part will be filled in by IJGT
\begin{center}
Communicated by\;
\end{center}
%------------------------------------------------------------------------------------%

\begin{abstract} 

We summarize several results about non-simplicity, solvability and normal structure of finite groups related to the number of conjugacy classes appearing in the product or the power of conjugacy classes. We also collect some problems that have only been partially solved. 
\end{abstract}

\vskip 0.2 true cm

%------------------------------------------------------------------------------------%

\pagestyle{myheadings}
\markboth{\rightline {\sl Int. J. Group Theory x no. x (201x) xx-xx \hskip 7 cm  Beltr\'an, Felipe and Melchor }}
         {\leftline{\sl Int. J. Group Theory x no. x (201x) xx-xx \hskip 7 cm  Beltr\'an, Felipe and Melchor}}

\bigskip
%------------------------------------------------------------------------------------%
%------------------------------------------------------------------------------------%
\section{Introduction}
Let $G$ be a finite group. A classic problem about conjugacy classes is the study of the group structure taking into account the information provided for the product of classes (see for instance \cite{products}). Particularly, since the product of conjugacy classes of $G$ is a $G$-invariant set, then it is a union of classes and accordingly, there exist several researches about the structure of a finite group regarding the number of conjugacy classes appearing in the product of its classes. Some of these are related to the normal structure, the solvability or the non-simplicity of the group. The aim of this survey is to summarize some new results and several open problems in this research topic.\\

This overview is divided into four parts. In the second part we present the advances made on Arad and Herzog's conjecture during the past few years, which still remains open. It asserts that in a non-abelian simple group the product of two conjugacy classes cannot be a conjugacy class. In the third section we will focus our attention on the case when the power of a class is again a class. In Section 4, we will address the problem of the power of a conjugacy class when it is a union of the trivial class and another class. Finally, in Section 5, we will outline the advances on the case in which the power of a class is a union of two classes being one of them the inverse of the other. Throughout this survey, we will approach to some characterizations in terms of irreducible characters.

\section{Products of conjugacy classes: Arad and Herzog's conjecture}
In \cite{products}, Z. Arad and M. Herzog conjectured one of the most significant problems concerning the product of two conjugacy classes.

\begin{problem}[Arad and Herzog's conjecture] 
In a non-abelian simple group the product of two conjugacy classes is not a conjugacy class.
\end{problem}

This was first checked for several families of simple groups such as for alternating groups, Suzuki groups, ${\rm PSL}_{2}(q)$, non-abelian simple groups of order less than one million, and 15 of the 26 sporadic simple groups (see \cite{products}). Since then, several authors have shown the conjecture to hold for certain groups. The following characterization in terms of irreducible characters can be used to check the conjecture, for example, for all the sporadic simple groups and for each group whose character table is known (for instance by using the Atlas \cite{Atlas} or {\sf GAP} \cite{GAP}).\\

We will denote by $x^G$ the conjugacy class of each element $x\in G$.

\begin{thm}\label{Moori}{\rm (Preliminaries of \cite{MooriViet})}
Let $G$ be a group and let $a$, $b$, $c$ $\in G$ be non-trivial elements of $G$. The following conditions are equivalent:
\begin{enumerate}
\item $a^Gb^G=c^G$
\item $\chi(a)\chi(b)=\chi(c)\chi(1)$ for all $\chi \in {\rm Irr}(G)$.
\end{enumerate}
\end{thm}

By using Theorem \ref{Moori}, J. Moori and H.P. Tong-Viet showed in \cite{MooriViet} that the conjecture is true for the following families of simple groups: ${\rm PSL}_{3}(q)$, ${\rm PSU}_{3}(q)$, where $q$ is a prime power,$^{2}{\rm G}_{2}(q)$ with $q=3^{2m+1}$, ${\rm PSp}_{4}(q)$ where $q$ is a prime-power,  ${\rm PSp}_{2n}(3)$ with $n\geq 2$, and  ${\rm PSU}_{n}(2)$, with $(n,3)=1$ and $n\geq 4$. More recently, R. Guralnick, G. Malle and P.H. Tiep validated Arad and Herzog's conjecture in several more cases \cite{GuralMallePham}. They also proved that if $G$ is a finite simple group of Lie type and $A$ and $B$ are non-trivial conjugacy classes, either both semisimple or both unipotent, then $AB$ is not a conjugacy class. In addition, they obtained a strong version of Arad and Herzog's conjecture for simple algebraic groups, and in particular, they show that almost always the product of two conjugacy classes in a simple algebraic group consists of infinitely many conjugacy classes. A weaker variation of Arad and Herzog's conjecture is posed in \cite{Nuestro7}. 

\begin{problem}[Conjecture 1.1 of \cite{Nuestro7}]\label{Con1}
In a non-abelian simple group the product of $n$ conjugacy classes, with $n$ a fixed natural number and $n\geq 2$, is not a conjugacy class.
\end{problem}

Let $K_1,\ldots, K_n$ be conjugacy classes with $n\geq 2$ such that the product $K_1K_2\cdots K_n$ is a conjugacy class $D$. Since the product of conjugacy classes is $G$-invariant, for each $i$ there exists a set of conjugacy classes $C_1,\ldots, C_t$ such that $K_1K_2\cdots K_{i-1}K_{i+1}\cdots K_n=C_1 \cup \ldots \cup C_t$ and we can write $K_1K_2\ldots K_n=K_iC_1\cup \cdots \cup K_iC_t$. Then $K_iC_j=D$ for all $i,j$ and hence, if Arad and Herzog's conjecture was proved, then we would automatically have Problem \ref{Con1} solved.\\

The authors provide a characterization in terms of characters of the situation in Problem \ref{Con1} that extends Theorem \ref{Moori}. In particular, such characterization holds for the case in which the power of a class is again a class.

\begin{thm}[Theorem 2.16 of \cite{Nuestro7}]\label{Char1} Let $K_1,\ldots, K_n $ be conjugacy classes of $G$ and write $K_i={x_i}^G$ with $x_i \in G$. Then $K_1K_2\cdots K_n=D$ where $D=d^G$ if and only if $$\chi(x_1)\cdots \chi(x_n)=\chi(1)^{n-1}\chi(d)$$ for all $\chi \in {\rm Irr}(G)$. In particular, if $K$ is a conjugacy class of $G$ and $x\in K$, then $K^n$ is a conjugacy class for some $n\in \mathbb{N}$ if and only if $$\chi(x)^n=\chi(1)^{n-1}\chi(x^n)$$ for all $\chi \in {\rm Irr}(G)$.
\end{thm}

Let $K$ be a class. When considering products of classes, a particular case is $KK^{-1}$ and the simplest situation is $KK^{-1}=1 \cup D$ where $D$ is a class. Under this hypothesis we obtain that the group is not simple by means of the Classification of the Finite Simple Groups (CFSG) and we conjecture that the subgroup $\langle K\rangle$ is solvable. Unfortunately, we have only been able to get such solvability in some particular cases.

\begin{thm} {\rm (Theorem C of \cite{Nuestro6})}\label{N1} Let $K$ be a conjugacy class of a finite group $G$ and suppose that $KK^{-1}= 1 \cup D$, where $D$ is a conjugacy class of $G$. Then $G$ is not simple, $|D|$ divides $|K|(|K|-1)$ and $\langle K\rangle/\langle D\rangle$ is cyclic. In addition,
\begin{enumerate}
\item If $|D|=|K|-1$, then $\langle K\rangle$ is metabelian. More precisely, $\langle D\rangle$ is $p$-elementary abelian for some prime $p$.
\item If $|D|=|K|$, then $\langle K\rangle$ is solvable with derived length at most 3.
\item If $|D|=|K|(|K|-1)$, then $\langle K\rangle$ is abelian.
\end{enumerate}
\end{thm}

We remark that $KK^{-1}=1 \cup D$ forces $D$ to be real, but it does not necessarily imply that $K$ is real too. Moreover, under the assumption of Theorem \ref{N1}, if $K$ is real, then it can be easily proved that $D$ is real too, and hence, $KK^{-1}=K^{2}=1\cup D$. This case will be addressed in Section 4.\\

The next natural and simplest case when considering the product of a class and its inverse has the form $KK^{-1}=1 \cup D \cup D^{-1}$. Also in this case it can be checked that the group is not simple. To prove this as well as to obtain the non-simplicity in Theorem \ref{N1}, the following characterization in terms of characters is used. 

\begin{thm}[Theorem B of \cite{Nuestro6}] Let $G$ be a group and $x, d\in G$. Let $K=x^G$ and $D=d^G$. The following assertions are equivalent:
\begin{enumerate}[label=\alph*)]
\item $KK^{-1}=1\cup D \cup D^{-1}$.
\item For every $\chi \in$ {\rm Irr}$(G)$ $$|K||\chi(x)|^2=\chi(1)^2+\frac{(|K|-1)}{2}\chi(1)(\chi(d)+\chi(d^{-1})).$$
\end{enumerate}

In particular, if $D=D^{-1}$, then $KK^{-1}=1\cup D$ if and only if for every $\chi \in$ {\rm Irr}$(G)$ $$|K||\chi(x)|^2=\chi(1)^2+(|K|-1)\chi(1)\chi(d).$$
\end{thm}

Under the hypotheses of Theorem \ref{N1} the group $G$ need
not be solvable. The typical non-solvable situation is a group of type $Z.S.2$, where $|Z| = 3$, $Z$ is in the center of $Z.S$, and in addition, $S$ is a non-solvable group acted by an automorphism of order 2 such that the non-trivial elements of $Z$ are conjugate by this automorphism.

\section{Power of class which is a class}
In \cite{GuralnickNavarro}, Guralnick and Navarro confirmed Arad and Hergoz's conjecture for the particular case of the square of a conjugacy class. They demonstrated that when the square of a conjugacy class of $G$ is again a class, then $G$ is not a non-abelian simple group. However, the result goes much further.

\begin{thm}[Theorem 2.2 of \cite{GuralnickNavarro}]\label{Gab1}
Let $G$ be a finite group and $K=x^G$ the conjugacy class of $x$ in $G$. The following assertions are equivalent:
\begin{enumerate}
\item $K^{2}$ is a conjugacy class of $G$.
\item $K=x[x,G]$ and {\rm \textbf{C}}$_{G}(x)=$ {\rm \textbf{C}}$_{G}(x^2)$.
\item $\chi(x)=0$ or $|\chi(x)|=\chi(1)$ for all $\chi\in$ {\rm Irr}$(G)$, and {\rm \textbf{C}}$_{G}(x)=$ {\rm \textbf{C}}$_{G}(x^2)$.
\end{enumerate}
In this case, $[x,G]$ is solvable, and $\langle K \rangle= \langle x\rangle [x, G]$ too. Furthermore, if the order of $x$ is a power of a prime $p$, then $[x,G]$ has normal $p$-complement.
\end{thm}

Observe that $K^2=K$ can never happen since, by Theorem 3.1, this implies that $K=x[x, G]=K^{2}=x^2[x, G]$. This means that $x\in [x, G]$ and $K=[x,G]$, a contradiction. Furthermore, there is an equivalent property to the fact that the square of a conjugacy class is again a class in terms of characters. That is, $\chi(x)^2=\chi(1)\chi(x^2)$ for every $\chi \in {\rm Irr}(G)$ if and only if $K^2$ is a conjugacy class where $K=x^G$, with $x\in G$. This is a particular case of Theorem \ref{Moori}.\\

The most relevant fact of Theorem \ref{Gab1} is that the authors obtained the solvability of the normal subgroup $[x, G]$ (which coincides with $\langle K\rangle$) by means of the CFSG. However, the assertion about the case of a class of $p$-elements does not require the CFSG. The key fact for proving the solvability is that all elements of $x[x, G]$ are $G$-conjugate. Moreover, the following result is used. 

\begin{thm} [Theorem 3.2 of \cite{GuralnickNavarro}]\label{key} Let $G$ be a finite group and let $N$ be a normal subgroup of $G$. Let $x\in G$ such that all elements of $xN$ are conjugate in $G$. Then
\begin{enumerate}
\item $N$ is solvable.
\item If $\pi$ is the set of prime divisors of $o(x)$, then $x$ normalizes some Hall $\pi$-complement $H$ of $N$ on which it acts fixed point freely.
\item If $x$ is a $p$-element for some prime $p$, then $N$ has normal $p$-complement.  
\end{enumerate}
\end{thm}

With regard to Problem \ref{Con1}, the assertion is proved to be true for the particular case of the power of a conjugacy class. We find a normal subgroup of a group when there is a conjugacy class satisfying that some of its powers is again a conjugacy class, and obtain an equivalent property in terms of the irreducible characters of the group. This result is in fact a generalization of Theorem \ref{Gab1} and its proof is similar.

\begin{thm}[Theorem 2.2 of \cite{Nuestro7}] \label{1} Let $K=x^G$ with $x\in G$, $n\in \mathbb{N}$ and $n\geq 2$. The following assertions are equivalent:
\begin{enumerate}[label=(\alph*)]
\item $K^n$ is a conjugacy class.
\item ${\rm \bf C}_G(x)={\rm \bf C}_G(x^n)$ and $N=x^{-1}K=K^{-1}K=[x,G]\unlhd G$.
\item ${\rm \bf C}_G(x)={\rm \bf C}_G(x^n)$ and $\chi(x)=0$ or $|\chi(x)|=\chi(1)$ for all $\chi\in$ {\rm Irr}$(G)$.
\end{enumerate}
\end{thm}

Notice that under the hypotheses of the above theorem, $G$ is not a non-abelian simple group. This happens because we know that $K=xN$ and if $N=1$, then $x$ is central, and if $N=G$, then $K=G$, so $N$ is always a non-trivial proper subgroup of $G$.\\

The following corollaries are obtained as a consequence of the above theorem. 

\begin{cor}[Corollary 2.4 of \cite{Nuestro7}]\label{C1} Let $K=x^G$ with $x \in G$ such that $K^n$ is a conjugacy class for some $n\in \mathbb{N}$ with $n\geq 2$. Then $|K^r| = |K|$ for all $r\in \mathbb{N}$ and $K^{o(x)+1}=K$ and $K^{o(x)-1}=K^{-1}$. Moreover, $K^m$ is a conjugacy class for all $m\in \mathbb{N}$ such that $(m, o(x))=1$.
\end{cor}

\begin{cor}[See Corollary 2.5 and Theorem 2.12 of \cite{Nuestro7}]\label{C2} If $K=x^G$ with $x\in G$ such that $K^n=K$ for some $n\in \mathbb{N}$ with $n\geq 2$, then:

\begin{enumerate}[label=(\alph*)]
\item $K^{k(n-1)+r}=K^r$ for every $r, k\in \mathbb{N}$. 
\item $K^{n-1}=[x, G]\unlhd G$. 
\item $\pi(o(x))\subseteq \pi(n-1)$, where $\pi(t)$ denotes the set of primes dividing the number $t$.
\item If $K$ is real, then $x$ is a $2$-element and $K^m=K$ for every odd number $m$. Also, $K^2=[x,G]\unlhd G$.
\end{enumerate}
\end{cor}

\begin{cor}[Corollary 2.8 of \cite{Nuestro7}]\label{C3} Let $G$ be a finite group and let $\pi$ be a set of primes. Suppose that for each conjugacy class $K$ of $\pi$-elements of $G$ there exists $n \in \mathbb{N}$ such that $K^n$ is a conjugacy class. Then $G/${\rm \textbf{O}}$_{\pi'}(G)$ is nilpotent. In particular, if $\pi=\pi(G)$, then $G$ is nilpotent.
\end{cor}

The authors also obtain the solvability of the subgroup generated by a class when one of its powers is a class by employing similar arguments to those of Guralnick and Navarro, which appeal to the CFSG. Theorem \ref{1} joint with Theorem \ref{key} are the main ingredients to prove the next theorem.

\begin{thm}[Theorem 2.15 of \cite{Nuestro7}] \label{A} Let $K=x^G$ be a conjugacy class of a group $G$. If there exists $n\in \mathbb{N}$ and $n\geq 2$ satisfying that $K^n$ is a conjugacy class, then $\langle K\rangle$ is solvable.
\end{thm}

The authors make emphasis in the cases where the CFSG is not necessary to obtain the solvability in Theorem \ref{A}. This occurs for classes of elements of prime order, classes of $2$-elements and real classes.

\begin{thm}[Theorem 2.11 of \cite{Nuestro7}]\label{NOCFSG} Let $K$ be a conjugacy class of an element $x \in G$. Suppose that there exists $n\in \mathbb{N}$ with $n\geq 2$ satisfying that $K^n$ is a conjugacy class. Then:
\begin{enumerate}
\item If $o(x)$ is a prime $p$, then $\langle K\rangle$ is solvable.
\item If $x$ is a $2$-element, then $\langle K\rangle$ is solvable. 
\item If $K^n=D$ where $D$ is a real conjugacy class, then $\langle K\rangle$ is solvable. Also, $D^3=D$ and $D$ is a class of a $2$-element.
\end{enumerate}
\end{thm}

In the particular case of Theorem \ref{NOCFSG}(3), if in addition $n=2^a$ for some $a\in \mathbb{N}$, then more information can be given. 

\begin{prop}[See Theorem 2.12 of \cite{Nuestro7}]
 Let $G$ be a finite group and let $K$ be a conjugacy class of an element $x\in G$. If $K^n=D$ where $D$ is a real conjugacy class with $n=2^a$ for some $a\in \mathbb{N}$, then $|K|$ is odd and $o(x)=2^{a+1}$.
\end{prop}

Observe that if a real conjugacy class $K$ satisfies that there exists $n\in \mathbb{N}$ such that $K^n=D$ where $D$ is a conjugacy class, then $D$ is also a real class. However, if a class $K$ satisfies that there exists $n\in \mathbb{N}$ such that $K^n=D$ with $D$ real, then $K$ need not be real. A trivial example of this situation occurs in $\mathbb{Z}_4$.

\section{Power of a class which is a union of the trivial class and another class}

As we have explained in Section 2, if $K$ and $D$ are conjugacy classes of a group $G$ such that $KK^{-1}=1 \cup D$, then $G$ is not simple (see \cite{Andrews1} for more details on this problem). A natural question is what happens if the power of a class is union of the trivial class and another one, so we establish the following problem. 

\begin{problem}\label{Con2} Let $G$ be a finite group and let $K$ be a conjugacy class. If $K^n=1 \cup D$ for some $n\in \mathbb{N}$ and $n\geq 2$, and some conjugacy class $D$, then $\langle K \rangle$ is solvable. In particular, $G$ is not a non-abelian simple group.  
\end{problem}

It can be easily proved that if $K^n=1 \cup D$ for some $n \in \mathbb{N}$, then $KK^{-1}=1 \cup D$, so by Theorem \ref{N1} the group is not simple. However, the converse is not true. Let $G=SL(2,3)$ and let $K$ be one of the two conjugacy classes of elements of order 6, which satisfies $|K|=4$. It follows that $KK^{-1}=1 \cup D$ where $D$ is the unique conjugacy class of size $6$. However, there is no $n\in \mathbb{N}$ with $K^n=1 \cup D$.\\

Problem \ref{Con2} is proved in \cite{Nuestro5} for the particular case $n=2$ without using the CFSG, and the structure of $\langle K\rangle$ and $\langle D\rangle$ are completely determined. Observe that if $K^n=1 \cup D$ for some $n \in \mathbb{N}$ and $K$ is real, then $K^{2}=1 \cup D$.

\begin{thm}[Theorem A of \cite{Nuestro5}]\label{A1} Let $K=x^G$ be a conjugacy class of a finite group $G$ and suppose that $K^2= 1 \cup D$, where $D$ is a conjugacy class of $G$. Then $\langle D\rangle=[x,G]$ is either cyclic or a $p$-group for some prime $p$, and so $\langle K\rangle=\langle x\rangle[x,G]$ is solvable. More precisely,

\begin{enumerate}
\item Suppose that $|K|=2$.
\begin{enumerate}
\item If $o(x)=2$, then $\langle K\rangle \cong \mathbb{Z}_{2}\times \mathbb{Z}_{2}$ and $\mathbb{Z}_{2}\cong \langle D\rangle\subseteq$ {\rm \textbf{Z}}$(G)$. 
\item If $o(x)=n>2$, then $\langle K\rangle \cong \mathbb{Z}_{n}$ and $\langle D\rangle$ is cyclic.
\end{enumerate}
\item Suppose that $|K|\geq 3$.
\begin{enumerate}
\item If $o(x)=2$, then either $\langle K \rangle$ and  $\langle D\rangle$ are 2-elementary abelian groups or $\langle D\rangle$ is a $p$-group and $|K|=p^r$ with $p$ an odd prime and $r$ a positive integer.
\item If $o(x)>2$, then $\langle D\rangle$ is a $p$-elementary abelian group for some odd prime $p$. Furthermore, either $o(x)=p$ or $o(x)=2p$. 
\end{enumerate}
 In every case, $|\langle K\rangle/\langle D\rangle|\leq 2$.
\end{enumerate}
\end{thm}

All cases of Theorem \ref{A1} are possible and examples of each of them are given in \cite{Nuestro5}. The techniques for proving Theorem \ref{A1} are relatively elementary although Glauberman's ${\rm Z}^*$ theorem \cite{Glauberman} and a result of Y. Berkovich and L. Kazarin in \cite{BerkoKaza} are used. Both require tools from modular representation theory, so Theorem \ref{A1} is based on it as well. Other two main ingredients of the proof of Theorem \ref{A1} are Burnside's classification of finite 2-groups having exactly one involution and the classification of groups of order 16.\\

The following theorem shows that Problem \ref{Con2} is completely solved. 

\begin{thm}[Theorem B of \cite{Nuestro7}]\label{B} Let $K=x^G$ with $x\in G\setminus \{1\}$ such that $K^n=1 \cup D$ for some $n\in \mathbb{N}$ and $n\geq 2$ where $D$ is a conjugacy class of $G$. Then $KK^{-1}=1 \cup D$ and $\langle K\rangle$ is solvable.
\end{thm}

To look for examples we use {\sf GAP} \cite{GAP} and particularly the {\sc SmalGroups} library \cite{SmallGroups}. The $m$-th group of order $n$ in the {\sc SmallGroups} library is identified by $n\#m$.

\begin{exmp}
Let us show two examples of the case $K^n=1\cup D$ with $n=3$. In the first, we have $x^n=1$ and in the second $x^n\neq 1$. Let $G=A_4$ and $K=(1\, 2\, 3)^G$, which satisfies $|K|=4$ and $o((1\, 2\, 3))=3$. Furthermore, $K^3=1 \cup D$ where $D=((1\, 2)(3\, 4))^G$. On the other hand, let $G=(\mathbb{Z}_7\rtimes \mathbb{Z}_9)\rtimes \mathbb{Z}_2=126\#11$ having a conjugacy class $K$ of elements of order 21 satisfying that $K^3=1\cup D$ and $|K|=6$ where $D$ is a class of elements of order 7 and $|D|=6$. In this example, $\langle K\rangle\cong \mathbb{Z}_{21}$.
\end{exmp}

The authors also provide a characterization with irreducible characters of the property stated in Problem \ref{Con2}, which is the following.

\begin{thm}[Theorem 3.9 of \cite{Nuestro7}]\label{Char2} Let $G$ be a finite group and let $K=x^G$. Then $K^n=1 \cup D$ where $D=d^G$ if and only if there exist positive integers $m_1$ and $m_2$ such that $$\chi(x)^n|K|^n=\chi(1)^{n-1}(m_1\chi(1)+m_2|D|\chi(d))$$ for all $\chi \in {\rm Irr}(G)$ and $|K|^n=m_1+m_2|D|$, where

$$m_1=\frac{|K|^n}{|G|}\sum_{\chi \in {\rm Irr}(G)}\frac{\chi(x)^n}{\chi^{n-2}(1)} \, \, \, \, \, \,  and \, \, \, \, \, \, m_2=\frac{|K|^n}{|G|}\sum_{\chi \in {\rm Irr}(G)}\frac{\chi(x)^n\overline{\chi(d)}}{\chi^{n-1}(1)}.$$
\end{thm}

\section{Power of a class which is a union of a class and its inverse}
A variant of the problem addressed in Section 4 is to study what happens when the power of a conjugacy class is a union of two classes, one of them being the inverse of the other. It is believed the following to hold in this case.

\begin{problem}[Conjecture 1.2 of \cite{Nuestro7}]\label{Con3} Let $G$ be a finite group and let $K$ be a conjugacy class. If $K^n=D \cup D^{-1}$ for some $n\in \mathbb{N}$ with $n\geq 2$ and $D$ a conjugacy class, then $\langle K \rangle$ is solvable. In particular, $G$ is not simple.  
\end{problem}

Observe that if $K=x^G$ with $x\in G$ and $K^n=D \cup D^{-1}$ for some $n \in \mathbb{N}$ with $D \neq D^{-1}$, then $K$ is non-real. Indeed, suppose that $K$ is real and $x^n \in D$. We have $x^{-1}=x^g$ for some $g \in G$. Then $(x^n)^g=(x^g)^n=(x^{-1})^n=x^{-n} \in D\cap D^{-1}$, so $D=D^{-1}$, that is, $D$ is real.\\

The authors provide the following evidence to support the conjecture of Problem \ref{Con3}.

\begin{thm}[Theorem C of \cite{Nuestro7}] \label{C} Let $G$ be a finite group and let $K$ be a conjugacy class. If $K^n=D \cup D^{-1}$ for some $n\in \mathbb{N}$ and $n\geq 2$ and $D$ a conjugacy class, then either $|D|=|K|/2$ or $|K|=|D|$. In the first case, $\langle K\rangle$ is solvable.
\end{thm}

\begin{exmp}
We show that both cases of Theorem \ref{C} are possible. Let $G=\mathbb{Z}_8 \rtimes \mathbb{Z}_2=M_{16}=\langle a, x\, \, | \, \, a^8=x^2=1, \, \, a^{x^{-1}}=a^5\rangle$ and $K=a^G$. We have $K^2= D \cup D^{-1}$, $|K|=2$ and $|D|=1$. On the other hand, let $G=\mathbb{Z}_2 \times (\mathbb{Z}_7 \rtimes \mathbb{Z}_3)=42\#2$ and $K=x^G$ where $o(x)=14$. We have $K^2= D\cup D^{-1}$ and $|K|=|D|=3$.
\end{exmp}

The conjecture of Problem \ref{Con3} is also demonstrated when $D=K$ by simply working with properties of the complex group algebra $\mathbb{C}[G]$. The following result, which is a particular case of the product of two classes for which Arad and Herzog's conjecture holds, is also used to prove Theorem \ref{D}.

\begin{lem}[Lemma 3.1 of \cite{Nuestro7}]\label{L1}
Let $G$ be a group and $K$, $L$ and $D$ are non-trivial conjugacy classes of $G$ such that $KL=D$ with $|D|=|K|$. Then $G$ possesses a solvable proper normal group which is $\langle LL^{-1}\rangle$. 
In particular, $\langle L\rangle$ is solvable.
\end{lem}

\begin{thm}[Theorem D of \cite{Nuestro7}]\label{D} Let $G$ be a finite group and let $K=x^G$ be a conjugacy class of $G$. If $K^2=K \cup K^{-1}$, then $\langle K\rangle$ is solvable. Moreover, $x$ is a $p$-element for some prime $p$.
\end{thm}

The following property is useful to check that the Conjecture in Problem \ref{Con3} is true for certain groups from their character tables. In particular, this happens for the sporadic simple groups.

\begin{thm}\label{Char3} Let $G$ be a finite group and let $K$ be a conjugacy class of an element $x\in G$. Then $K^n=D \cup D^{-1}$ where $D$ is a conjugacy class if and only if there exist positive integers $m_1$ and $m_2$ such that $$\chi(x)^n|K|^n=\chi(1)^{n-1}|D|(m_1\chi(x^n)+m_2\chi(x^{-n}))$$ for all $\chi \in {\rm Irr}(G)$ and $|K|^n=(m_1+m_2)|D|$ where

$$m_1=\frac{|K|^n}{|G|}\sum_{\chi \in {\rm Irr}(G)}\frac{\chi(x)^n\overline{\chi(x^n)}}{\chi^{n-1}(1)} \, \, \, \, \, \, and \, \, \, \, \, \, m_2=\frac{|K|^n}{|G|}\sum_{\chi \in {\rm Irr}(G)}\frac{\chi(x)^n\chi(x^{n})}{\chi^{n-1}(1)}.$$

In particular, $$\chi(x)^n+\chi(x^{-1})^n=\chi(1)^{n-1}(\chi(x^n)+\chi(x^{-n}))$$ for all $\chi \in {\rm Irr}(G)$.
\end{thm}

\begin{rmk} Recall that the smallest integer $m$ satisfying $C^m=G$ for each non-trivial conjugacy class $C$ of $G$ is called the covering number of $G$. The covering number of each sporadic simple group is at most 6 {\rm (}\cite{Zisser} and \cite{products}{\rm )}. It can be checked by using the character tables {\rm (}included in {\sf GAP}{\rm )} that for any of these groups, any two non-trivial conjugacy classes of it and $n<6$, there exists an irreducible character that does not satisfy the equation of Theorem \ref{Char3}.
\end{rmk}

\section*{Acknowledgments}
Some results of this survey were obtained during the stay of C. Melchor at the University of Cambridge in autumn 2017, which was financially supported by the grant E-2017-02, Universitat Jaume I of Castell\'on. She would like to thank R. Camina and the Department of Mathematics for their warm hospitality. A. Beltr\'an and M.J. Felipe are supported by the Valencian Government, Proyecto PROMETEOII/2015/011.

%-----------------------------------------------------------------------------
%-----------------------------------------------------------------------------

\bigskip
\bigskip

{\footnotesize \pn{\bf Antonio Beltr\'an}\; \\ {Departamento de Matem\'aticas}, {Universidad Jaume I, 12071,} {Castell\'on, Spain}\\
{\tt Email: abeltran@uji.es}\\

{\footnotesize \pn{\bf Mar\'{\i}a Jos\'e Felipe}\; \\ {Instituto Universitario de Matem\'atica Pura y Aplicada}, { Universitat Polit\`ecnica de Val\`encia, 46022,} {Valencia, Spain}\\
{\tt Email: mfelipe@mat.upv.es }\\

{\footnotesize \pn{\bf Carmen Melchor}\; \\ {Departamento de Matem\'aticas}, {Universidad Jaume I, 12071,} {Castell\'on, Spain}\\
{\tt Email: cmelchor@uji.es}\\
\end{document}